\documentclass[12pt, journal, draftclsnofoot, onecolumn]{IEEEtran}
\ifCLASSOPTIONcompsoc
\else
\fi

\ifCLASSINFOpdf
\else
\fi
\usepackage[parfill]{parskip}    
\usepackage{amsmath, amssymb, latexsym, enumerate, mathrsfs}
\usepackage{graphicx}
\usepackage{hyperref}
\usepackage{graphicx}
\usepackage{array}
\usepackage{times}

\usepackage{moreverb}
\usepackage{marvosym}
\usepackage{color,soul}
\definecolor{lightblue}{rgb}{.90,.95,1}
\sethlcolor{yellow}

\usepackage{hyperref}
\DeclareMathOperator*{\arginf}{arg\,inf}

\newtheorem{lemma}{\bf Lemma}
\newtheorem{proposition}{\bf Proposition}

\newtheorem{remark}{Remark}

\newtheorem{assumption}{\bf Assumption}
\newtheorem{definition}{\bf Definition}

\def\etal{\mbox{et al.}}

\begin{document}
%
\title{On the stochastic decision problems with backward stochastic viability property}

\author{Getachew K.~Befekadu,~\IEEEmembership{Member,~IEEE} 
\IEEEcompsocitemizethanks{\IEEEcompsocthanksitem G. K. Befekadu is with the National Research Council, 
Air Force Research Laboratory \& Department of Industrial System Engineering, University of Florida - REEF, 1350 N. Poquito Rd, Shalimar, FL 32579, USA. \protect\\
E-mail: gbefekadu@ufl.edu}}

\IEEEcompsoctitleabstractindextext{%
\begin{abstract}
In this paper, we consider a stochastic decision problem for a system governed by a stochastic differential equation, in which an optimal decision is made in such a way to minimize a vector-valued accumulated cost over a finite-time horizon that is associated with the solution of a certain multi-dimensional backward stochastic differential equation (BSDE). Here, we also assume that the solution for such a multi-dimensional BSDE {\it almost surely} satisfies a backward stochastic viability property w.r.t. a given closed convex set. Moreover, under suitable conditions, we establish the existence of an optimal solution, in the sense of viscosity solutions, to the associated system of semilinear parabolic PDEs. Finally, we briefly comment on the implication of our results.

\end{abstract}

\begin{IEEEkeywords}
Diffusion processes, forward-backward SDEs, stochastic decision problem, value functions, viscosity solutions.
\end{IEEEkeywords}}

\maketitle

\IEEEdisplaynotcompsoctitleabstractindextext

%
\IEEEpeerreviewmaketitle

\section{Introduction} \label{S1}
Let $\bigl(\Omega, \mathcal{F}, \mathbb{P}, \{\mathcal{F}_t \}_{t \ge 0}\bigr)$ be a probability space, and let $\{B_t\}_{t \ge 0}$ be a $d$-dimensional standard Brownian motion, whose natural filtration, augmented by all $\mathbb{P}$-null sets, is denoted by $\{\mathcal{F}_t\}_{t \ge 0}$, so that it satisfies the {\it usual hypotheses} (e.g., see \cite{Pro90}). We consider the following system governed by a stochastic differential equation (SDE)
\begin{align}
d X_t = f\bigl(t, X_t, u_t\bigr) dt + \sigma\bigl(t, X_t, u_t\bigr)dB_t, \quad X_0=x, \quad  0 \le t \le T,  \label{Eq1.1} 
\end{align}
where
\begin{itemize}
\item $X_{\cdot}$ is an $\mathbb{R}^{d}$-valued diffusion process,
\item $u_{\cdot}$ is a $U$-valued measurable admissible control process (where $U$ is a compact set in $\mathbb{R}^{d}$) such that for all $t > s$, $(B_t - B_s)$ is independent of $u_r$ for $r \le s$ (nonanticipativity condition) and
\begin{align*}
\mathbb{E} \int_{s}^{t} \vert u_t\vert^2 dt < \infty, \quad \forall t \ge s,
\end{align*}
\item the function $f \colon [0, \infty) \times \mathbb{R}^d \times U \rightarrow \mathbb{R}^{d}$ is uniformly Lipschitz, with bounded first derivative, and
\item $\sigma \colon [0, \infty) \times \mathbb{R}^{d} \times U \rightarrow \mathbb{R}^{d \times d}$ is Lipschitz with the least eigenvalue of $\sigma\,\sigma^T$ uniformly 
bounded away from zero for all $(x,u) \in \mathbb{R}^{d} \times U$ and $t \ge 0$.
\end{itemize}
In this paper, we specifically consider a stochastic decision problem for the above system, where the admissible decision $u_{\cdot}$ is a $U$-valued measurable control process from the set $\mathcal{U}_{[0,T]}$ with
\begin{align}
\mathcal{U}_{[0,T]} \triangleq \Bigl\{u \colon [0,T] \times \Omega \rightarrow U \,\Bigl\vert  u \,\, \text{is an} \,\, \bigl\{\mathcal{F}_t\bigr\}_{t\ge 0}\text{- adapted and}\,\, \mathbb{E} \int_{0}^{T} \vert u_t\vert^2 dt < \infty \Bigr\}. \label{Eq1.2} 
\end{align}
Furthermore, we consider the following vector-valued cost functional, which provides information on the accumulated cost on the time interval $[0, T]$,
\begin{align}
  \xi_{0,T}(u) = \int_0^T c\bigl(t, X_t, u_t\bigr) dt + \Psi(X_T), \label{Eq1.3} 
\end{align}
where the cost-rate $c \colon [0,T] \times \mathbb{R}^d \times U \rightarrow \mathbb{R}^n$ and the final-stage cost $\Psi \colon \mathbb{R}^d \rightarrow \mathbb{R}^n$ are measurable functions. 

Here, we remark that the corresponding solution $X_t=X_t^{0,x;u}$ depends on the admissible decision $u_{\cdot} \in \mathcal{U}_{[0,T]}$ as well as on the initial condition $X_0=x$.  As a result of this, for any time-interval $[t, T]$, with $t \in [0, T]$, the vector-valued accumulated cost $\xi_{t,T}$ depends on the admissible decision $u_{\cdot} \in \mathcal{U}_{[t,T]}$. Moreover, we also assume that $f$, $\sigma$, $c$ and $\Psi$ satisfy the following condition 
\begin{align}
\bigl\vert f\bigl(t, x, u) \bigr\vert + \bigl\vert \sigma\bigl(t, x, u \bigr) \bigr\vert + \bigl\vert c\bigl(t, x, u\bigr) \bigr\vert + \bigl\vert \Psi\bigl(x\bigr) \bigr\vert \le C \bigl(1 + \bigl\vert x \bigr\vert^p + \bigl\vert u \bigr\vert \bigr) \label{Eq1.4}
\end{align}
for all $\bigl(t, x, u \bigr) \in [0,T] \times \mathbb{R}^{d} \times U$, $p \ge 1$ and for some constant $C > 0$.

Next, let us introduce the following measurable spaces that will be useful throughout the paper. For any Euclidean space $H$, we denote by $L_{ad}^2\bigl(\Omega, C([0, T], H)\bigr)$ the closed linear subspace of adapted processes of $L^2\bigl(\Omega, \mathcal{F}_t, \mathbb{P}, C([0,T], H)\bigr)$, and $L_{ad}^2\bigl(\Omega \times [0, T], H\bigr)$ is the Hilber space of adapted measurable stochastic processes $\varphi \colon \Omega \times ]0, T[ \rightarrow H$ such that $\bigl\Vert \varphi \bigr\Vert_{[t, T]}^2 = \mathbb{E}\bigl\{ \int_t^T \bigl\vert \varphi_s \bigr \vert^2 ds  \bigr\}< \infty$.

On the same probability space $\bigl(\Omega, \mathcal{F}, \mathbb{P}, \{\mathcal{F}_t \}_{t \ge 0})$, we consider the following multi-dimensional BSDE 
\begin{align}
- d Y_t = G\bigl(t, Y_t, Z_t\bigr) dt - Z_tdB_t, \quad Y_T=\xi, \label{Eq1.5}
\end{align}
where the terminal value $Y_T=\xi$ belongs to $L^2\bigl(\Omega, \mathcal{F}_T, \mathbb{P}, \mathbb{R}^n\bigr)$ and the function $G \colon \Omega \times [0, T] \times \mathbb{R}^n \times \mathscr{L}(\mathbb{R}^{d}, \mathbb{R}^{n})\rightarrow \mathbb{R}^n$, with property that $G\bigl(\cdot, \cdot, y, z\bigr)$ is progressively measurable and $G\bigl(\Omega, \cdot, y, z\bigr)$ is continuous. Moreover, we also assume the following conditions on $G$.

\begin{assumption} \label{AS1} ~\\\vspace{-0.15in}
\begin{enumerate} [(i)]
\item $G$ is Lipschitz in $(y, z)$, i.e., there exists a constant $C > 0$, such that
\begin{align*}
\bigl\vert G\bigl(t, y, z\bigr) - G\bigl(t, y', z'\bigr) \bigr\vert \le C \bigl(\bigl\vert y - y' \bigr\vert + \bigl\Vert z - z' \bigr\Vert\bigr).
\end{align*}
\item $\sup_{t \le T} G\bigl(t, 0, 0\bigr) \in L^2\bigl(\Omega, \mathcal{F}, \mathbb{P} \bigr)$.
\end{enumerate}
for all $(t, y, z, y', z')$, $\mathbb{P}$-{\it almost everywhere} on $\Omega$.
\end{assumption}

\begin{remark} \label{RR1}
In this paper, we also require all constituent functions $G_j(t, y, z)$ for $j=1,2, \ldots, n$ (with $G(t, y, z)= [G_1(t, y, z), G_2(t, y, z), \dots, G_n(t, y, z)]^T$) to have the following special  structure $G_j(t, y, z) =  G_j(t, y, z^j)$ and $z^j$ is the $j$th row of the $z \in \mathbb{R}^{n \times d}$ for $j=1,2, \ldots, n$ (cf. the proof part of Proposition~\ref{P3}).
\end{remark}

Next, we state the following lemma, which is used to establish the existence of unique adapted solutions (e.g., see \cite{ParP90} for additional discussions).
\begin{lemma} \label{L1}
Suppose that Assumption~\ref{AS1} holds. Then, for any $\xi \in L^2\bigl(\Omega, \mathcal{F}_T, \mathbb{P}, \mathbb{R}^n\bigr)$, the BSDE in \eqref{Eq1.5}, with terminal condition $Y_T=\xi$, i.e.,
\begin{align}
 Y_t = \xi + \int_t^T G\bigl(s, Y_s, Z_s\bigr) ds - \int_t^TZ_s dB_s, \quad 0 \le t \le T  \label{Eq1.6}
\end{align}
has a unique adapted solution
\begin{align}
 \bigl(Y_t^{T,g,\xi}, Z_t^{T,g,\xi}\bigr)_{0 \le t \le T} \in L_{ad}^2\bigl(\Omega, C([0, T], \mathbb{R}^n) \bigr) \times L_{ad}^2\bigl(\Omega \times ]0, T[, \mathscr{L}(\mathbb{R}^d, \mathbb{R}^n) \bigr). \label{Eq1.7}
\end{align}
\end{lemma}

In what follows, we introduce a definition for conditional $G$-expectation, which is associated with the progressively measurable function $G$ of BSDE in \eqref{Eq1.5}.
\begin{definition} \label{Df1} 
For any $\xi \in L^2\bigl(\Omega, \mathcal{F}_T, \mathbb{P}, \mathbb{R}^n\bigr)$, let 
\begin{align*}
 \bigl(Y_t^{T,g,\xi}, Z_t^{T,g,\xi}\bigr)_{0 \le t \le T} \in L_{ad}^2\bigl(\Omega, C([0, T], \mathbb{R}^n) \bigr) \times L_{ad}^2\bigl(\Omega \times ]0, T[, \mathscr{L}(\mathbb{R}^d, \mathbb{R}^n) \bigr)
\end{align*}
be the unique solution to the BSDE in \eqref{Eq1.5} with terminal condition $Y_T=\xi  \in L^2\bigl(\Omega, \mathcal{F}_T, \mathbb{P}, \mathbb{R}^n\bigr)$. Then, we define the conditional $G$-expectation of $\xi$ as follows
\begin{align}
\mathcal{E}^G \bigl[\xi \bigl \vert \mathcal{F}_t \bigr] \triangleq Y_t^{T,g,\xi}, \quad t \in [0, T].  \label{Eq1.8}
\end{align}
\end{definition}
 
Note that such a nonlinear expectation (i.e., the conditional $G$-expectation) is widely used for evaluating the risk of uncertain future outcomes or costs, where the classical average performance criteria (which is based on the standard linear expectation) may not be sufficient to account for how risks are perceived by decision makers.

Here, it is worth mentioning that some interesting studies on risk measures, based on the conditional $G$-expectation, have been reported in the literature (e,g. see \cite{Pen04}, \cite{CorHMP02} and \cite{Ros06} for establishing connections between the risk measures and the function $G$ of BSDE). Moreover, such risk measures are widely used for evaluating the risk of uncertain future outcomes or costs, and also assisting with stipulating minimum interventions for risk management (e.g., see \cite{ArtDEH99}, \cite{Pen04}, \cite{DetS05} or \cite{CorHMP02} for related discussions). On the other hand, for decision problems involving multi-dimensional BSDEs coupled with forward-SDEs, there are some interesting studies, based on the stochastic backward viability property, that establish condition for the solutions of the associated PDEs (e.g., see \cite{Xu16} and \cite{BuckQR00} for additional discussions; see also \cite{AuDaPra95} and \cite{AuDaPra98} on the notion of viability properties for SDEs and inclusions). Note that the rationale behind our framework, which follows in some sense the settings of these papers, is to show how a backward stochastic viability property can be systematically used to obtain consistently optimal decision solutions.

The remainder of this paper is organized as follows. In Section~\ref{S2}, using the basic remarks made in Sections~\ref{S1}, we state the problem of optimal decisions for the system governed by a (forward) stochastic differential equation. In Section~\ref{S3}, we present our main result -- where we establish the existence of optimal decisions, in the sense of viscosity solutions, to the associated system of semilinear parabolic partial differential equations. Finally, Section~\ref{S4} provides further remarks.

\section{Problem formulation} \label{S2}
In order to make our problem formulation mathematically more appealing, for any $(t, x) \in [0, T] \times \mathbb{R}^d$, we consider the following forward SDE with an initial condition $X_t^{t,x;u} = x$
\begin{align}
d X_s^{t,x;u} = f\bigl(t, X_s^{t,x; u}, u_s\bigr) dt  + \sigma\bigl(t, X_s^{t,x; u}, u_s\bigr)dB_t, \quad t \le s \le T, &  \label{Eq2.1}
\end{align}
where $u_{\cdot}$ is a $U$-valued measurable control process from the set $\mathcal{U}_{[t,T]}$.

Let $\xi$ be a real-valued random variable from $L^2(\Omega, \mathcal{F}_T, \mathbb{P}, \mathbb{R}^n)$ and we further suppose that the data $\xi$ takes the following form
\begin{align}
 \xi = \Psi(X_T^{t,x;u_{\cdot}}),\quad \mathbb{P}-{almost \, surely \, (a.s)}.  \label{Eq2.2}
 \end{align}
Moreover, we introduce the following value function\footnote{$V^{u}\bigl(t, x\bigr)$ is a vector-valued function, i.e., 
\begin{align*}
 V^{u}\bigl(t, x\bigr) =& [V_1^{u}\bigl(t, x\bigr), V_2^{u}\bigl(t, x\bigr), \ldots, V_n^{u}\bigl(t, x\bigr)]^T,  \,\, (t, x) \in [0, T] \times \mathbb{R}^d \,\, \&  \,\, u \in U.
 \end{align*}}
\begin{align}
 & \quad  V^{u}\bigl(t, x\bigr) = \mathcal{E}^{G} \bigl[\xi_{t,T} \bigl(u\bigr) \bigl \vert \mathcal{F}_t \bigr] \quad  \text{with} \quad   u_{\cdot} \in \mathcal{U}_{[t,T]}, \label{Eq2.3}
\end{align}
where
\begin{align}
\xi_{t,T}\bigl(u\bigr) = \int_t^T c\bigl(s, X_s^{t,x;u}, u_s\bigr) ds + \Psi(X_T^{t,x;u}). \label{Eq2.4}
\end{align}
Note that we can express the above value function as follow
\begin{align}
 V^{u}\bigl(t, x\bigr) &= \xi_{t,T}\bigl(u\bigr) + \int_t^T G\bigl(s, Y_s^{t,x;u}, Z_s^{t,x;u}\bigr) ds - \int_t^T Z_s^{t,x;u} dB_s \notag \\
                                                                      &= \Psi(X_T^{t,x;u}) + \int_t^T \Bigl\{c\bigl(s, X_s^{t,x;u}, u_s\bigr)  + G\bigl(s, Y_s^{t,x;u}, Z_s^{t,x;u} \bigr) \Bigr\}ds - \int_t^T Z_s^{t,x;u}dB_s, \label{Eq2.5}
\end{align}
where the function $G$ is assumed to satisfy Assumption~\ref{AS1}. Furthermore, noting the condition in \eqref{Eq1.4}, then $\bigl(Y_s^{t,x;u}, Z_s^{t,x;u}\bigr)_{t \le s \le T}$ is an adapted solution on $[t, T] \times \Omega$ and belongs to $L_{ad}^2\bigl(\Omega, C([0, T], \mathbb{R}^n) \bigr) \times L_{ad}^2\bigl(\Omega \times ]0, T[, \mathscr{L}(\mathbb{R}^d, \mathbb{R}^n) \bigr)$. Equivalently, we can rewrite \eqref{Eq2.5} as a multi-dimensional BSDE on the probability space $\bigl(\Omega, \mathcal{F}, \mathbb{P}, \{\mathcal{F}_t \}_{t \ge 0})$, i.e.,
\begin{align}
 d Y_s^{t,x;u} = \hat{G}\bigl(s, X_s^{t,x;u}, Y_s^{t,x;u}, Z_s^{t,x;u}\bigr) ds - Z_s^{t,x;u} dB_s, \quad s \in [t, T], \quad Y_T^{t,x;u} = \Psi(X_T^{t,x;u}). \label{Eq2.6}
\end{align}
where 
\begin{align*}
 \hat{G}\bigl(s, X_s^{t,x;u}, Y_s^{t,x;u}, Z_s^{t,x;u}\bigr) = c\bigl(s, X_s^{t,x;u}, u_s\bigr)  + G\bigl(s, Y_s^{t,x;u}, Z_s^{t,x;u}\bigr).
\end{align*}
with $\hat{G}_j\bigl(s, x, y, z\bigr)=\hat{G}_j\bigl(s, x, y, z^j\bigr)$ for $j=1,2, \ldots, n$.

Note that the problem of finding an optimal decision $u_{\cdot}^{\ast} \in \mathcal{U}_{[t,T]}$, with $t \in [0, T]$, that minimizes the vector-valued accumulated cost is then reduced to finding an optimal solution for
\begin{align}
\inf_{u_{\cdot} \in \mathcal{U}_{[t,T]}}  J\bigr[u\bigl], \label{Eq2.7}
\end{align}
where 
\begin{align}
J\bigr[u\bigl] = \mathcal{E}^{G} \bigl[\xi_{t,T}\bigl(u\bigr) \bigl \vert \mathcal{F}_t \bigr] . \label{Eq2.8}
\end{align}
Here, we remark that, for any given $u_{\cdot} \in \mathcal{U}_{[t,T]} $, with $t \in [0, T]$, if the forward-backward stochastic differential equations (FBSDEs) in \eqref{Eq2.1} and \eqref{Eq2.6} admit weak solutions, then the corresponding solution $X_s^{t,x;u}$ depends uniformly on $u_{\cdot} \in \mathcal{U}_{[t,T]}$, for $s \in [t, T]$.

Let $K$ be a closed convex set in $\mathbb{R}^n$, then we recall the notion of viability property for the BSDE in \eqref{Eq2.6}.
\begin{definition} \label{Df2}
Let  $u_{\cdot} \in \mathcal{U}_{[0,T]}$ be an admissible decision process, then, for a nonempty closed convex set $K \subset \mathbb{R}^n$, we have
\begin{enumerate} [(a)]
\item A stochastic process $\bigl\{Y_t^{0,x;u_{\cdot}}, \,\,t \in [0, T] \bigr\}$ is viable in $K$ if and only if for $\mathbb{P}$-{\it almost}\, $\omega \in \Omega$ 
\begin{align}
Y_t^{0,x;u_{\cdot}}(\omega) \in K, \quad \forall t \in [0, T]. \label{Eq2.9}
\end{align}
\item The closed convex set $K$ enjoys the backward stochastic viability property (BSVP) for the equation in \eqref{Eq2.6} if and only if for all $\tau \in [0, T]$, with equation~\eqref{Eq2.2}, i.e., 
\begin{align}
\forall \, \xi \in L^2\bigl(\Omega, \mathcal{F}_{\tau}, \mathbb{P}, \mathbb{R}^n\bigr), \label{Eq2.10}
\end{align}
there exists a solution pair $\bigl(Y_{\cdot}^{0,x;u_{\cdot}}, Z_{\cdot}^{0,x;u_{\cdot}} \bigr)$ to the BSDE in \eqref{Eq2.6} over the time interval $[0, \tau]$,
\begin{align*}
Y_s^{0,x;u_{\cdot}} = \xi + & \int_s^{\tau} \hat{G}\bigl(r, X_r^{0,x;u_{\cdot}}, Y_r^{0,x;u_{\cdot}}, Z_r^{0,x;u_{\cdot}}\bigr) dr - \int_s^{\tau} Z_r^{0,x; u} dB_r, 
\end{align*}
with 
\begin{align*}
 \bigl(Y_{\cdot}^{0,x;u_{\cdot}}, Z_{\cdot}^{0,x;u_{\cdot}} \bigr) \in  L_{ad}^2\bigl(\Omega, C([0, \tau], \mathbb{R}^n) \bigr) \times L_{ad}^2\bigl(\Omega \times ]0, \tau[, \mathscr{L}(\mathbb{R}^d, \mathbb{R}^n) \bigr)
\end{align*}
such that $\bigl\{Y_s^{0,x;u_{\cdot}}, \,\,s \in [0, \tau] \bigr\}$ is viable in $K$.
\end{enumerate}
\end{definition}

With respect to the above convex set $K$, let us define the projection of a point $a$ onto $K$ as follow
\begin{align}
\Pi_{K} (a) = \Bigl\{b \in K \, \bigl\vert \, \vert a - b \vert = \min _{c \in K} \vert a - c \vert = d_{K}(a) \Bigr\}. \label{Eq2.12}
\end{align}
Notice that, since $K$ is convex, from the Motzkin's theorem, $\Pi_{K}$ is single-valued. Further, we recall that $d_{K}^2(\cdot)$ is convex; and thus, due to Alexandrov's theorem \cite{Ale39}, $d_{K}^2(\cdot)$ is almost everywhere twice differentiable.

Moreover, on the space $C_b^{1,2}([t, T] \times \mathbb{R}^d; \mathbb{R}^n)$, for any $(t,x) \in [0, T] \times \mathbb{R}^d$, with $\varphi \in C_b^{1,2}([t, T] \times \mathbb{R}^d; \mathbb{R}^n)$, and, noting the statement in Remark~\ref{RR1}, we consider the following system of semilinear parabolic partial differential equations (PDEs)
\begin{eqnarray}
\left.\begin{array}{r}
 \dfrac{\partial \varphi_j(t, x)}{\partial t}  + \inf_{u \in U} \Bigl\{\mathcal{L}_{t}^{u} \varphi_j(t, x)  + \hat{G}_j\bigl(t, \varphi(t, x), D_x \varphi_j(t, x) \cdot \sigma(t, x, u)\bigr)\Bigr\} = 0 \\
    j =1,2, \ldots, n
\end{array}\right\}  \label{Eq2.13}
\end{eqnarray}
with the following boundary condition
\begin{align}
\varphi(T, x) &= \Psi(x), \quad x \in \mathbb{R}^d,  \label{Eq2.14}
\end{align}
where, for any $\phi(x) \in C_0^{\infty}(\mathbb{R}^d)$, the second-order linear operators $\mathcal{L}_{t}^{u}$ are given by
\begin{align}
 \mathcal{L}_{t}^{u} \phi(x) = \dfrac{1}{2} \operatorname{tr} \Bigl\{a(t, x, u) D_{x}^2 \phi(x)\Bigr\} &+ f(t, x, u) D_{x} \phi(x), \,\,  t \in [0, T], \label{Eq2.15}
\end{align}
with $a(t, x, u) = \sigma(t, x, u) \sigma^T(t, x, u)$, $D_{x}$ and $D_{x}^2$, (with $D_{x}^2 = \bigl({\partial^2 }/{\partial x_k \partial x_l} \bigr)$) are the gradient and the Hessian (w.r.t. the variable $x$), respectively.

Here, we remark that the above system of semilinear parabolic PDEs in \eqref{Eq2.13} together with the boundary condition of \eqref{Eq2.14}, is associated with the stochastic decision problem in \eqref{Eq2.7} (see also equations~\eqref{Eq2.17} and \eqref{Eq2.18} below), restricted to $\Sigma_{[t,T]}$ (cf. Definition~\ref{Df4}). Note that the problem of FBSDEs and the solvability of the associated system of semilinear parabolic PDEs have been well studied in literature (e.g., see \cite{HUP95}, \cite{LIW14},  \cite{ParT99}, \cite{MaZZ08}, \cite{MaPY94} and \cite{Pen92}).

Next, let us define the viability property for the system of semilinear parabolic PDEs in \eqref{Eq2.13} as follow.
\begin{definition} \label{Df4}
The system of semilinear parabolic PDEs in \eqref{Eq2.13} enjoys the viability property w.r.t. the closed convex set $K$ if and only if, for any $\Psi \in C_p(\mathbb{R}^d; \mathbb{R}^n)$ taking values in $K$, the viscosity solution to \eqref{Eq2.13} satisfies
\begin{align}
\forall (t, x) \in [0, T] \times \mathbb{R}^d, \quad \varphi(t,x) \in K. \label{Eq2.16}
\end{align}
\end{definition}
Further, we introduce the following definition for an admissible decision system $\Sigma_{[t, T]}$, for $t \in [0, T]$, which provides a logical construct for our main results (e.g., see also \cite{LIW14}).
\begin{definition} \label{Df5}
For a given finite-time horizon $T>0$, we call $\Sigma_{[t, T]}$, with $t \in [0, T]$, an admissible decision system, if it satisfies the following conditions:
\begin{itemize}
\item $\bigl(\Omega, \mathcal{F},\{\mathcal{F}_t \}_{t \ge 0}, \mathbb{P}\bigr)$ is a complete probability space.
\item $\bigl\{B_s\bigr\}_{s \ge t}$ is a $d$-dimensional standard Brownian motion defined on $\bigl(\Omega, \mathcal{F}, \mathbb{P}\bigr)$ over $[t, T]$ and $\mathcal{F}^t \triangleq \bigl\{\mathcal{F}_s^t\bigr\}_{s \in [t, T]}$, where $\mathcal{F}_s^t = \sigma\bigl\{\bigl(B_s; \,t \le s \le T \bigr)\bigr\}$ is augmented by all $\mathbb{P}$-null sets in $\mathcal{F}$.
\item $u_{\cdot} \colon \Omega \times [s, T]  \rightarrow U$ is an $\bigl\{\mathcal{F}_s^t\bigr\}_{s \ge t}$-adapted process on $\bigl(\Omega, \mathcal{F}, \mathbb{P}\bigr)$ with 
\begin{align*}
\mathbb{E} \int_{s}^{T} \vert u_{\tau} \vert^2 d \tau < \infty,  \quad s \in [t, T].
\end{align*}
\item For any $x \in \mathbb{R}^d$, the FBSDEs in \eqref{Eq2.1} and \eqref{Eq2.6} admit a unique solution set \, $ \bigl\{X_{\cdot}^{s,x;u_{\cdot}}, Y_{\cdot}^{s,x;u_{\cdot}}, Z_{\cdot}^{s,x;u_{\cdot}}\bigr\}$ on $\bigl(\Omega, \mathcal{F}, \mathcal{F}^t, \mathbb{P}\bigr)$ and $Y_{\cdot}^{s,x;u_{\cdot}}(\omega) \in K$ for $\mathbb{P}$-{\it almost} $\omega \in \Omega$ for all $s \in [t, T]$.
\end{itemize}
\end{definition}

Then, with restriction to the above admissible decision system $\Sigma_{[0, T]}$, we can state the optimal decision problem as follows.

{\it Problem}: Find an optimal decision $u_{\cdot}^{\ast} \in \mathcal{U}_{[0,T]}$ such that
\begin{align}
  u_{\cdot}^{\ast} \in \arginf J\bigr[ u \bigl] \Bigl \vert \,u_{\cdot} \in \mathcal{U}_{[0,T]} \,\text{restricted to} \,\Sigma_{[0, T]}  \label{Eq2.17}
\end{align}
Furthermore, the optimal vector-valued accumulated cost $J$ over the time-interval $[0, T]$ is given  
\begin{align}
J\bigr[u^{\ast}\bigl] &= \mathcal{E}^{G} \left [ \int_0^T c\bigl(s, X_s^{0,x;u^{\ast}}, u_s\bigr) ds + \Psi(X_T^{0,x;u^{\ast}})\Bigl \vert \mathcal{F}_0, \right] \notag \\
                               &\equiv  \int_0^T c\bigl(s, X_s^{0,x;u^{\ast}}, u_s\bigr) ds + \Psi(X_T^{0,x;u^{\ast}}),\,\, X_0^{0,x;u^{\ast}} = x. \label{Eq2.18}
\end{align}

In the following section, assuming the Markovian framework, we establish the existence of an optimal solution, in the sense of viscosity solutions (e.g., see \cite{CraIL92} or \cite{FleS06} for additional discussions on the notion of viscosity solutions), for the above stochastic decision problem with restriction to $\Sigma_{[0, T]}$. 

\section{Main results} \label{S3}
In this section, we present our main results -- where we make use of the following observations: If the closed convex set $K \in \mathbb{R}^n$ enjoys the BSVP for the multi-dimensional BSDE of \eqref{Eq2.6}. Then, it is sufficient for the existence of an optimal decision for the stochastic decision problem in \eqref{Eq2.17}; provided that the viscosity solutions  for the corresponding system of semilinear parabolic PDEs in \eqref{Eq2.13} together with \eqref{Eq2.14} enjoy the viability property w.r.t. the same closed convex set $K$.

Let us first state the following two propositions, i.e., Proposition~\ref{P1} and Proposition~\ref{P2} (whose proofs for one dimensional BSDEs are also given in \cite{BefVP16}), that will be useful for proving our main results.
\begin{proposition} \label{P1}
Suppose Assumption~\ref{AS1} together with \eqref{Eq1.4} hold. Then, for any $(t,x) \in [0, T] \times \mathbb{R}^d$ and $u_{\cdot} \in \mathcal{U}_{[t,T]}$, the system equations (i.e., the FBSDEs) in \eqref{Eq2.1} and \eqref{Eq2.6} admit unique adapted solutions
\begin{eqnarray}
\left.\begin{array}{c}
X_{\cdot}^{t,x;u} \in L_{ad}^2\bigl(\Omega, C([t, T], \mathbb{R}^n) \bigr)\\
\bigl(Y_{\cdot}^{t,x;u}, Z_{\cdot}^{t,x;u}\bigr) \in L_{ad}^2\bigl(\Omega, C([t, T], \mathbb{R}^n) \bigr) \times L_{ad}^2\bigl(\Omega \times ]t, T[, \mathscr{L}(\mathbb{R}^d, \mathbb{R}^n) \bigr)
\end{array}\right\}  \label{Eq3.1}
\end{eqnarray}
Furthermore, the value function $V^{u}\bigl(t, x\bigr)$ is deterministic.
\end{proposition}

\begin{IEEEproof} [{\it Proof of Proposition~\ref{P1}}]
Notice that $f$ and $\sigma$ are bounded and Lipschitz continuous w.r.t. $(t,x) \in [0, T] \times \mathbb{R}^d$ and uniformly for $u \in U$. Then, for any $(t,x) \in [0, T] \times \mathbb{R}^d$ and progressively measurable process $u_{\cdot}$, there always exists a unique path-wise solution $X_{\cdot}^{t,x;u_{\cdot}} \in L_{ad}^2\bigl(\Omega, C([t, T], \mathbb{R}^n) \bigr)$ for the forward SDE in \eqref{Eq2.1}. On the other hand, consider the following multi-dimensional BSDE,
\begin{align}
 -d \hat{Y}_s^{t,x;u_{\cdot}} = \hat{G}\bigl(s, X_s^{t,x;u_{\cdot}}, \hat{Y}_s^{t,x;u_{\cdot}}, Z_s^{t,x;u_{\cdot}}\bigr) ds - Z_s^{t,x;u_{\cdot}} dB_s, \label{EqA1.1}
\end{align}
where 
\begin{align*}
\hat{Y}_T^{t,x;u_{\cdot}} = \int_t^T c\bigl(\tau, X_{\tau}^{t,x;u_{\cdot}}, u_{\tau}\bigr) d\tau + \Psi(X_T^{t,x;u_{\cdot}}). 
\end{align*}
From Lemma~\ref{L1}, equation~\eqref{EqA1.1} admits unique solutions $\bigl(\hat{Y}_{\cdot}^{t,x;u_{\cdot}}, Z_{\cdot}^{t,x;u_{\cdot}}\bigr)$ in $\in L_{ad}^2\bigl(\Omega, C([t, T], \mathbb{R}^n) \bigr) \times L_{ad}^2\bigl(\Omega \times ]t, T[, \mathscr{L}(\mathbb{R}^d, \mathbb{R}^n) \bigr)$. Furthermore, if we introduce the following
\begin{align*}
Y_s^{t,x;u_{\cdot}} = \hat{Y}_s^{t,x;u_{\cdot}} - \int_t^s c\bigl(\tau, X_{\tau}^{t,x;u_{\cdot}}, u_{\tau}\bigr) d\tau,  \quad s \in [t, T].
\end{align*}
Then, the forward version of the BSDE in \eqref{Eq2.6} holds with $\bigl(Y_{\cdot}^{t,x;u_{\cdot}}, Z_{\cdot}^{t,x;u_{\cdot}}\bigr)$. Moreover, we also observe that $Y_t^{t,x;u_{\cdot}}$ is deterministic. This completes the proof of Proposition~\ref{P1}.
\end{IEEEproof}

\begin{proposition} \label{P2}
Let  $(t,x) \in [0, T] \times \mathbb{R}^d$ and $u_{\cdot} \in \mathcal{U}_{[t,T]}$ be restricted to $\Sigma_{[t, T]}$. Then, for any $r \in [t, T]$ and $\mathbb{R}^d$-valued $\mathcal{F}_r^t$-measurable random variable $\eta$, we have
\begin{align}
 V^{u}\bigl(r, \eta\bigr) =  \mathcal{E}^G\Bigl[\int_r^T c\bigl(s, X_s^{r,\eta;u}, u_s\bigr) ds + \Psi(X_T^{r,\eta;u}) \Bigl \vert \mathcal{F}_r \Bigr],\quad P{\text-a.s.}  \label{Eq3.2}
\end{align}
\end{proposition}

\begin{IEEEproof} [{\it Proof of Proposition~\ref{P2}}]
For any $r \in [t, T]$, with $t \in [0, T]$, we consider the following probability space $\bigl(\Omega, \mathcal{F}, \mathbb{P}\bigl(\cdot\vert \mathcal{F}_r^t\bigr), \{\mathcal{F}^t\}\bigr)$ and notice that $\eta$ is deterministic under this probability space. Then, for any $s \ge r$, there exist progressively measurable process $\psi$ such that
\begin{align}
 u_s(\Omega)&= \psi(\Omega, B_{\cdot \wedge s}(\Omega)),\notag \\  
                        &= \psi(s, \bar{B}_{\cdot \wedge s}(\Omega) + B_{r}(\Omega)), \label{EqA2.1}
\end{align}
where $\bar{B}_s = B_s - B_r$ is a standard $d$-dimensional brownian motion. Note that $u_{\cdot}$is an $\mathcal{F}_r^t$-adapted process, then we have the following restriction w.r.t. $\Sigma_{[t, T]}$
\begin{align}
 \bigl(\Omega, \mathcal{F}, \{\mathcal{F}^t\}, \mathbb{P}\bigl(\cdot \vert \mathcal{F}_r^t\bigr)(\omega'), B_{\cdot}, u_{\cdot}\bigr) \in \Sigma_{[t, T]}, \label{EqA2.2}
\end{align}
where $\omega' \in \Omega'$ such that $\Omega' \in \mathcal{F}$, with $\mathbb{P}(\Omega')=1$. Furthermore, noting Lemma~\ref{L1}, if we work under the probability space $\bigl(\Omega', \mathcal{F}, \mathbb{P}\bigl(\cdot\vert \mathcal{F}_r^t\bigr)\bigr)$, then the statement in \eqref{Eq3.2} holds $\mathbb{P}$-{\it almost surely}. This completes the proof of Proposition~\ref{P2}. 
\end{IEEEproof}

\begin{proposition} \label{P3}
Let $u_{\cdot} \in \mathcal{U}_{[t,T]}$ be restricted to $\Sigma_{[t, T]}$, with $t \in [0, T]$. Suppose that the system of semilinear parabolic PDEs in \eqref{Eq2.13} enjoys the viability property w.r.t. the closed convex set $K$. Then, there exists a constant $C > 0$ such that $d_{K}^2(\cdot)$ is twice differentiable at $y$ and
\begin{align}
\bigl\langle y - \Pi_{K}(y),\, \hat{G}(t, x, y, z\sigma(t,x, u^{\neg j})) \bigr\rangle & \le \frac{1}{4}\bigl\langle D^2(d_{K}^2(y)) z \cdot \sigma(t,x, u),\,z \cdot \sigma(t,x, u) \bigr\rangle + C d_{K}^2(y), \notag \\
& \quad\quad \quad \forall (t, x, y, z) \in [0, T] \times \mathbb{R}^d \times \mathbb{R}^n \times \mathcal{L}(\mathbb{R}^d; \mathbb{R}^n). \label{Eq3.3} 
\end{align}
\end{proposition}

\begin{IEEEproof} [{\it Proof of Proposition~\ref{P3}}] 
The proof for the above proposition (which is an adaptation of \cite{BuckQR00}) involves a standard approximation procedure for the multi-dimensional BSDE in \eqref{Eq2.6}, with
\begin{align*}
Z_s^{t,x;u_{\cdot}} \in & \operatorname{span} \Bigl\{z \cdot \sigma(t, X_s^{t,x;u_{\cdot}}, u)  \, \bigl \vert \, z \in \mathcal{L}(\mathbb{R}^d; \mathbb{R}^n)\Bigr\}, \quad ds \otimes d \mathbb{P}-a.e. \,\, \text{on} \,\, [t, T],  \,\, 0 \le t \le T,
\end{align*}
and a further requirement for the closed convex set $K$ to enjoy the BSVP for the equation in \eqref{Eq2.6} (i.e., the adapted solution $\bigl\{Y_s^{t,x;u}, \, s \in [t, T] \bigr\}$ to be viable in $K$).
Here, we can show that the adapted solution $\bigl\{Y_s^{t,x;u_{\cdot}}, \, s \in [t, T] \bigr\}$ is viable in $K$, when the statement in \eqref{Eq3.3} holds true. Note that $K$ is a nonempty closed convex subset in $\mathbb{R}^n$. From \cite{DelZo95}, if $d_{K}^2$ is twice differentiable almost everywhere, then $\Pi_{K}$ is a single-valued mapping and 
\begin{eqnarray}
\left\{\begin{array}{l}
 \triangledown d_{K}^2(y) = 2(y - \Pi_{K}(y)), \quad \forall y \in  \mathbb{R}^n \\
 \Vert \Pi_{K}(y) - \Pi_{K}(x+y) \Vert \le \Vert x \Vert, \quad \forall (x,y) \in  \mathbb{R}^n \times \mathbb{R}^n.
 \end{array}\right. \label{EqA3.1}
\end{eqnarray}
Moreover, let $\Lambda_{K}$ be the set of all points of $\mathbb{R}^n$, where $d_{K}^2$ is twice differentiable.\footnote{Notice that this set is a Lebesque measure.} Further, let the measurable mapping $D^2(d_{K}^2) \colon \Lambda_{K} \rightarrow S(\mathbb{R}^n)$ be defined by the second order development of  $d_{K}^2$ in $y \in \Lambda_{K}$, i.e.,
\begin{align}
d_{K}^2(x+y) = d_{K}^2(y) + \langle \triangledown d_{K}^2(y),\,x\rangle + \frac{1}{2} \langle D^2(d_{K}^2(y))x,\,x\rangle + \alpha(y,x), \label{EqA3.2} 
\end{align}
where
\begin{align*}
 \frac{1}{\Vert x \Vert^2} \alpha(y, x)  \rightarrow  0 \quad \text{as} \quad x \rightarrow 0.
\end{align*}
Here, we claim that
\begin{eqnarray}
\left\{\begin{array}{l}
 (i) \quad 0 \le \frac{1}{2} D^2(d_{K}^2(y)) \le I_{n \times n}, \quad \forall y \in \Lambda_{K}, \\
 (ii) \quad \vert \alpha(y, x) \vert \le \Vert x \Vert^2, \quad \forall (x,y) \in \mathbb{R}^n \times \Lambda_{K}.
 \end{array}\right. \label{EqA3.3}
\end{eqnarray}
In order to verify the above conditions, let us first fix $y \in \Lambda_{K}$. Note that, since $d_{K}^2$ is convex, we have the following
\begin{align*}
d_{K}^2(x+y) - d_{K}^2(y) - \langle \triangledown d_{K}^2(y),\,x\rangle \ge 0, \quad \forall  x \in \mathbb{R}^n 
\end{align*}
On the other hand, we also have
\begin{align*}
d_{K}^2(x+y) - d_{K}^2(y) - \langle \triangledown d_{K}^2(y),\,x\rangle& \le \Vert (y +x) - \Pi_{K}(y) \Vert^2 - \Vert y  - \Pi_{K}(y) \Vert^2 -2 \langle x,\, y  - \Pi_{K}(y) \rangle \\
                                                                                                                &= \quad \Vert x \Vert^2. 
\end{align*}
Hence, from the above two inequalities, if we substitute $te$ for $x$, with $t > 0$, where $e$ is an arbitrary unit vector in $\mathbb{R}^n$. Then, we obtain the following
\begin{align*}
- \frac{1}{t^2} \langle \alpha(y, te) &\le \frac{1}{t^2} \bigl(d_{K}^2(y + te) - d_{K}^2(y) - t\langle \triangledown d_{K}^2(y),\,e\rangle \bigr) -  \frac{1}{t^2}\alpha(y, te)\\
                                                      &= \frac{1}{2} \langle D^2(d_{K}^2(y))e,\,e \rangle\\
                                                       &\le 1 - \frac{1}{t^2}\alpha(y, te). 
\end{align*}
Passing to the limit $t \rightarrow 0^{+}$, this gives the condition (i) of \eqref{EqA3.3}, which further implies (ii).

Let $\eta \in C^{\infty}(\mathbb{R}^n)$ be a nonnegative function with support in the unit ball and such that
\begin{align*}
\int_{\mathbb{R}^n} \eta(x) dx =1.
\end{align*}
For $\delta > 0$, we consider
\begin{align*}
\eta_{\delta}(x) dx =\frac{1}{\delta^n} \eta(\frac{1}{\delta} x)
\end{align*}
and
\begin{align*}
\phi_{\delta}(x) &= \int_{\mathbb{R}^n} d_{K}^2(x - \lambda) \eta_{\delta}(\lambda) d\lambda\\
           &\triangleq d_{K}^2(x) \star \eta_{\delta}(x), \quad x \in \mathbb{R}^n. 
\end{align*}
Notice that $\phi_{\delta} \in C^{\infty}(\mathbb{R}^n)$ and it also satisfies the following properties
\begin{eqnarray}
\left\{\begin{array}{l}
 (a) \quad 0 \le \phi_{\delta}(x) \le (d_{K}(x) + \delta)^2, \\
 (b) \quad  \triangledown \phi_{\delta}(x) = \int_{\mathbb{R}^n} \triangledown(d_{K}^2)(\lambda) \eta_{\delta}(x - \lambda) d\lambda,\\
      \quad \quad  \quad \Vert  \triangledown \phi_{\delta}(x) \Vert \le 2(d_{K}(x) + \delta),\\
 (c) \quad  D^2 \phi_{\delta}(x) = \int_{\mathbb{R}^n} D^2(d_{K}^2)(\lambda) \eta_{\delta}(x - \lambda) d\lambda,\\
       \quad \quad \quad 0 \le D^2 \phi_{\delta}(x) \le 2 I_{n \times n}, \quad \forall x \in \mathbb{R}^n.
 \end{array}\right. \label{EqA3.4}
\end{eqnarray}
Next, let us focus on properties (b) and (c) in \eqref{EqA3.4}, since property (a) is obvious. Note that $\Lambda_{K}$ is of full measure, then, for any $x$ and $\lambda$ belong to $\mathbb{R}^n$, we have
\begin{align*}
\phi_{\delta}(\lambda + x) - \phi_{\delta}(\lambda) & = \int_{\mathbb{R}^n} \bigl\{ d_{K}^2(y + x) - d_{K}^2(y) \bigr\} \eta_{\delta} (\lambda - y) dy \\
                                                                               & = \langle \int_{\mathbb{R}^n} \triangledown (d_{K}^2)(y) \eta_{\delta} (\lambda - y) dy, \, x \rangle + \frac{1}{2} \langle \bigl(\int_{\mathbb{R}^n} D^2 (d_{K}^2)(y) \eta_{\delta} (\lambda - y) dy\bigr)x, \, x \rangle + \varepsilon(\lambda, x),                                                                   
\end{align*}
where 
\begin{align*}
\varepsilon(\lambda, x) = \int_{\mathbb{R}^n} \alpha(y, x) \eta_{\delta} (\lambda - y) dy.
\end{align*}
Further, from Lebesque's dominance convergence theorem, we have the following
\begin{align*}
\frac{\varepsilon(\lambda, x)}{\Vert x \Vert^2} = \int_{\mathbb{R}^n} \frac{\alpha(y, x)}{\Vert x \Vert^2} \eta_{\delta} (\lambda - y) dy \rightarrow 0 \,\, \text{as} \,\, \Vert x \Vert \rightarrow 0, \,\, \forall \lambda \in \mathbb{R}^n.
\end{align*}
Next, consider $\xi \in L^2(\Omega, \mathcal{F}_T, \mathbb{P}; K)$ and let $\bigl(Y_{\cdot}^{0,x;u_{\cdot}}, Z_{\cdot}^{0,x;u_{\cdot}} \bigr)$ be a pair of unique adapted solutions to the following BSDE
\begin{align*}
Y_t^{0,x;u_{\cdot}} &= \xi + \int_t^T \hat{G}\bigl(s, X_s^{0,x;u_{\cdot}}, Y_s^{0,x;u_{\cdot}}, Z_s^{0,x;u_{\cdot}}\bigr) ds - \int_t^T Z_s^{0,x; u_{\cdot}} dB_s, \,\, t \in [0, T].
\end{align*}
From \eqref{EqA3.4}, we can apply the It\^{o} formula to $\phi_{\delta}(Y_t^{0,x;u_{\cdot}})$, for $0 \le t \le T$, with $\delta > 0$, then we obtain the following
\begin{align*}
 \mathbb{E} \phi_{\delta}(Y_t^{0,x;u_{\cdot}}) & = \mathbb{E} \xi +  \mathbb{E} \int_t^T \langle \triangledown \phi_{\delta}(Y_s^{0,x;u_{\cdot}}),\, \hat{G}\bigl(s, X_s^{0,x;u_{\cdot}} Y_s^{0,x;u_{\cdot}}, Z_s^{0,x;u_{\cdot}}\bigr) \rangle ds \\
& \quad  + \frac{1}{2}  \mathbb{E} \int_t^T \langle D^2(\phi_{\delta})(Y_s^{0,x;u_{\cdot}})Z_s^{0,x;u_{\cdot}},\, Z_s^{0,x;u_{\cdot}} \rangle ds \\
 & \,\, \le \delta^2 + \mathbb{E} \int_t^T \int_{\mathbb{R}^n} \Bigl\{ \langle \triangledown d_{K}^2(y),\, \hat{G}\bigl(s, X_s^{0,x;u_{\cdot}}, y, Z_s^{0,x;u_{\cdot}}\bigr) \rangle \\
&  \quad + \frac{1}{2}  \langle D^2(d_{K}^2(y))Z_s^{0,x;u_{\cdot}},\, Z_s^{0,x;u_{\cdot}} \rangle \Bigr\} \eta_{\delta} (Y_s^{0,x;u_{\cdot}} - y) dy ds \\
 & \quad - \mathbb{E} \int_t^T \int_{\mathbb{R}^n} \Bigl\{ \langle D^2 d_{K}^2(y),\, \hat{G}\bigl(s, X_s^{0,x;u_{\cdot}}, y, Z_s^{0,x;u_{\cdot}}\bigr)  \\
 & \quad\quad - \hat{G}\bigl(s, X_s^{0,x;u_{\cdot}}, Y_s^{0,x;u_{\cdot}}, Z_s^{0,x;u_{\cdot}}\bigr) \rangle \Bigr\}  \eta_{\delta} (Y_s^{0,x;u_{\cdot}} - y) dy ds.                                                                  
\end{align*}
Moreover, from \eqref{Eq3.3} and \eqref{EqA3.4}, for $\delta \in [0, T]$, we have the following
\begin{align*}
 \mathbb{E} \phi_{\delta}(Y_t^{0,x;u_{\cdot}}) & \le \delta^2 + C \mathbb{E} \int_t^T \int_{\mathbb{R}^n} d_{K}^2(y) \eta_{\delta} (Y_s^{0,x;u_{\cdot}} - y) dy ds \\
 & \quad\quad - \mathbb{E} \int_t^T \int_{\mathbb{R}^n} \Bigl\{ 2 d_{K}^2(y) \max_{y\colon \Vert Y_s^{0,x;u_{\cdot}} - y \Vert \le \delta}\Vert \hat{G}\bigl(s, X_s^{0,x;u_{\cdot}}, y, Z_s^{0,x;u_{\cdot}}\bigr) \\
 & \quad\quad\ -\hat{G}\bigl(s, X_s^{0,x;u_{\cdot}}, Y_s^{0,x;u_{\cdot}}, Z_s^{0,x;u_{\cdot}}\bigr) \Vert \Bigr\} \eta_{\delta} (Y_s^{0,x;u_{\cdot}} - y) dy ds \\                                                                 
 & \,\, \le \delta^2 + C \mathbb{E} \int_t^T  \mathbb{E} \phi_{\delta}(Y_t^{0,x;u_{\cdot}}) ds \\
 & \quad\quad - \mathbb{E} \int_t^T (1 + \phi_{\delta}(Y_t^{0,x;u_{\cdot}})) \max_{y\colon \Vert Y_s^{0,x;u_{\cdot}} - y \Vert \le \delta}\Vert \hat{G}\bigl(s, X_s^{0,x;u_{\cdot}}, y, Z_s^{0,x;u_{\cdot}}\bigr) \\
  & \quad \quad -\hat{G}\bigl(s, X_s^{0,x;u_{\cdot}}, Y_s^{0,x;u_{\cdot}}, Z_s^{0,x;u_{\cdot}}\bigr) \Vert ds.                                                                  
\end{align*}
Taking into account the function $g$ is uniformly continuous in its second variable, uniformly with respect to others, then we obtain that for some continuous increasing function $ \kappa: \mathbb{R}_{+} \rightarrow \mathbb{R}_{+}$, with $\kappa(0)=0$,
\begin{align}
 \mathbb{E} \phi_{\delta}(Y_t^{0,x;u_{\cdot}}) \le \delta^2 + \kappa(\delta) + (C+1) \mathbb{E} \int_t^T \phi_{\delta}(Y_s^{0,x;u_{\cdot}}) ds. \label{EqA3.5}
\end{align}
For $0 \le t \le T$, with small enough $\delta > 0$., from property (b) of \eqref{EqA3.4}, we obtain
\begin{align*}
 \mathbb{E} \phi_{\delta}(Y_t^{0,x;u_{\cdot}}) \le +\infty.
\end{align*}
This further allows us to apply Gronwall's inequality to \eqref{EqA3.5}. Hence, there is a real number $C > 0$, which does not depend on $t \in [0, T]$ and $\delta \in (0, 1)$, such that
\begin{align*}
 \mathbb{E} \phi_{\delta}(Y_t^{0,x;u_{\cdot}}) \le  C(\delta^2 + \kappa(\delta)), \quad 0 \le t \le T.
\end{align*}
Finally, since  $g$ is bounded, from Fatou's lemma and the dominated convergence theorem, we conclude that
\begin{align*}
 \mathbb{E} d_{K}^2(Y_t^{0,x;u_{\cdot}}) &\le \liminf_{\delta \rightarrow 0} \mathbb{E} \phi_{\delta}(Y_s^{0,x;u_{\cdot}}) \\
                                                                    &= 0, \quad 0 \le t \le T,
\end{align*}
i.e., $Y_t^{0,x;u_{\cdot}} \in K$, for any $t \in [0, T]$, $\mathbb{P}$-{\it almost everywhere}. This completes the proof of Proposition~\ref{P3}.
\end{IEEEproof}

\begin{remark}
Note that the above proposition provides an extension result for a multi-dimensional comparison theorem, based on the viability property w.r.t. the closed convex set $K$, for the solutions of the BSDE in \eqref{Eq2.6} on the probability space $\bigl(\Omega, \mathcal{F}, \mathbb{P}, \{\mathcal{F}_t \}_{t \ge 0})$ (e.g., see \cite{Xu16} for additional discussions)  
\end{remark}

In what follows, suppose that Proposition~\ref{P3} holds true, i.e., the system of semilinear parabolic PDEs in \eqref{Eq2.13} enjoys viability property w.r.t. the closed convex set $K$. Then, with restriction to $\Sigma_{[t, T]}$,  with $t \in [0, T]$, we can characterize the optimal decisions for the stochastic decision problem in \eqref{Eq2.17} as follows.
\begin{proposition} \label{P4}
Suppose that Proposition~\ref{P4} holds and let $\varphi \in C_b^{1,2}([0, T] \times \mathbb{R}^d; \mathbb{R}^n)$ satisfy \eqref{Eq2.13} with $\varphi\bigl(T, x\bigr)=\Psi(T, x)$ for $x \in \mathbb{R}^d$. Then, $\varphi\bigl(t, x\bigr) \le V^{u}\bigl(t, x\bigr)$ for any control $u_{\cdot} \in \mathcal{U}_{[t,T]}$ with restriction to $\Sigma_{[t, T]}$ and for all $(t,x) \in [0, T] \times \mathbb{R}^d$. Furthermore, if an admissible optimal decision process $u_{\cdot}^{\ast} \in \mathcal{U}_{[t,T]}$ exists, for almost all $(s, \omega) \in [0, T] \times \Omega$, together with the corresponding solution $X_s^{t,x; u^{\ast}}$, and satisfies
\begin{align}
 u_t^{\ast} \in\arginf_{u_{\cdot} \in \mathcal{U}_{[t,T]}}  \Bigl\{c_j\bigl(s, X_s^{0,x;u}, u_s\bigr) + & \mathcal{L}_{t}^{u} \varphi_j\bigl(s, X_s^{0,x;u}\bigr) +\hat{G}_j\bigl(s, D_x\varphi_j \bigl(s, X_s^{0,x;u}\bigr) \cdot \sigma\bigl(t, X_s^{0,x;u}, u_s\bigr) \bigr), \notag \\
 & \hspace{1.6in} j=1,2, \ldots, n \Bigr\}. \label{Eq3.15}
 \end{align}
Then, $\varphi\bigl(t, x\bigr) = V^{u^{\ast}}\bigl(t, x\bigr)$ for all $(t,x) \in [0, T] \times \mathbb{R}^d$.
\end{proposition}

\begin{IEEEproof} [{\it Proof of Proposition~\ref{P4}}] 
Suppose that $\varphi \in C_b^{1,2}([0, T] \times \mathbb{R}^d; \mathbb{R}^n)$  and, noting again the statement in Remark~\ref{RR1}, we assume that $\varphi_j \ge V_j^{u}$ on $[0, T] \times \mathbb{R}^d$ and $\max_{(t,x)} \bigl[V_j^{u}(t,x) - \varphi_j(t,x)\bigr] = 0$ for each $j=1,2, \ldots, n$. We consider a point $(t_{0},x_{0}) \in [0, T] \times \mathbb{R}^d$ so that $\varphi_j(t_{0},x_{0})= V_j^{u}(t_{0},x_{0})$ (i.e., a local maximum at $(t_{0},x_{0})$). Further, for a small $\delta t > 0$, we consider a constant control $u_s=\alpha$ for $s \in [t_{0},t_{0} + \delta t]$. Then, from \eqref{Eq3.2}, we have
\begin{align}
 \varphi_j(t_{0},x_{0}) = V_j^{u}(t_{0},x_{0}) & \le \mathcal{E}^{G_j}\Bigl[\int_{t_{0}}^{t_{0} + \delta t} c_j\bigl(s, X_s^{t_{0},x_{0};u}, \alpha \bigr) ds  + V_j^{w}(t_{0} + \delta t, X_{t_{0} + \delta t}^{t_{0},x_{0};u}) \Bigl \vert \mathcal{F}^{t_0}\Bigr]  \notag \\
                                       &\le \mathcal{E}^{G_j} \Bigl[\int_{t_{0}}^{t_{0} + \delta t} c_j\bigl(s, X_s^{t_{0},x_{0};u}, \alpha \bigr) ds  + \varphi_j(t_{0} + \delta t, X_{t_{0} + \delta t}^{t_{0},x_{0};u})\Bigl \vert \mathcal{F}^{t_0} \Bigr], \notag \\
                                       & \hspace{2.1in} j=1,2, \ldots, n. \label{EqP4.1}
\end{align}
Using the translation property of $\mathcal{E}^{G}[\,\cdot\,\vert \mathcal{F}^{t_0}]$, we obtain the following inequality
\begin{align}
\mathcal{E}^{G_j} \Bigl[\int_{t_{0}}^{t_{0} + \delta t} c_j\bigl(s, X_s^{t_{0},x_{0};u}, \alpha \bigr) ds + \varphi_j(t_{0} + \delta t, X_{t_{0} + \delta t}^{t_{0},x_{0};u}) - \varphi_j(t_{0},x_{0}) \Bigl \vert \mathcal{F}^{t_0} \Bigr] \ge 0. \label{EqP4.2}
\end{align}
Notice that $\varphi \in C_b^{1,2}([0, T] \times \mathbb{R}^d; \mathbb{R}^n)$, Then, using the It\^{o} formula, we can evaluate the difference between $\varphi_j(t_{0} + \delta t, X_{t_{0} + \delta t}^{t_{0},x_{0};u})$ and $\varphi_j(t_{0},x_{0})$, for $j=1,2, \ldots, n$ as follow
\begin{align}
\varphi_j(t_{0} + \delta t, X_{t_{0} + \delta t}^{t_{0},x_{0};u}) - \varphi_j(t_{0},x_{0}) & = \int_{t_{0}}^{t_{0} + \delta t} \Bigl[\dfrac{\partial}{\partial t} \varphi_j(s, X_{s}^{t_{0},x_{0};u}) + \mathcal{L}_{t}^{\alpha} \varphi_j(s, X_{s}^{t_{0},x_{0};u}) \Bigr] d s \notag  \\
  & \quad + \int_{t_{0}}^{t_{0} + \delta t} D_x \varphi_j(s, X_{s}^{t_{0},x_{0};u}) \cdot \sigma(s, X_{s}^{t_{0},x_{0};u}, \alpha) d B_s. \label{EqP4.3}
\end{align}
Moreover, if we substitute the above equation into \eqref{EqP4.2}, then we obtain the following
\begin{align}
 \mathcal{E}^{G_j} \Bigl[ \int_{t_{0}}^{t_{0} + \delta t} & \Bigl[c_j\bigl(s, X_s^{t_{0},x_{0};u}, \alpha \bigr) + \dfrac{\partial}{\partial t} \varphi_j(s, X_{s}^{t_{0},x_{0};u}) + \mathcal{L}_{t}^{\alpha} \varphi_j(s, X_{s}^{t_{0},x_{0};u}) \Bigr] d s \notag \\
& \, + \int_{t_{0}}^{t_{0} + \delta t} D_x \varphi_j(s, X_{s}^{t_{0},x_{0};u}) \cdot \sigma(s, X_{s}^{t_{0},x_{0};u}, \alpha) d B_s \Bigl \vert \mathcal{F}^{t_0} \Bigr]  \ge 0, \quad j=1,2, \ldots, n, \label{EqP4.4}
\end{align}
which amounts to solving the following system of BSDEs
\begin{align}
Y_{t_{0}}^{j;t_{0},x_{0};u} & = \int_{t_{0}}^{t_{0} + \delta t} \Bigl[c_j\bigl(s, X_s^{t_{0},x_{0};u}, \alpha \bigr) + \dfrac{\partial}{\partial t} \varphi_j(s, X_{s}^{t_{0},x_{0};u}) + \mathcal{L}_{t}^{\alpha} \varphi_j(s, X_{s}^{t_{0},x_{0};u}) \Bigr] d s \notag \\
& \quad + \int_{t_{0}}^{t_{0} + \delta t} D_x \varphi_j(s, X_{s}^{t_{0},x_{0};u}) \cdot \sigma(s, X_{s}^{t_{0},x_{0};u}, \alpha) d B_s \notag \\
&\quad \quad \quad + \int_{t_{0}}^{t_{0} + \delta t}\hat{G}_j\bigl(s, Z_s^{j;t_{0},x_{0};u}\bigr) ds - \int_{t_{0}}^{t_{0} + \delta t} Z_s^{j; t_{0},x_{0};u} dB_s, \quad j=1,2, \ldots, n.  \label{EqP4.5}
\end{align}
From Proposition~\ref{P1}, the above system of BSDEs admit unique solutions, i.e.,
\begin{align*}
Z_s^{j; t_{0},x_{0};u} & = D_x \varphi_j(s, X_{s}^{t_{0},x_{0};u}) \cdot \sigma(s, X_{s}^{t_{0},x_{0};u}, \alpha), \quad t_{0} \le s \le t_{0} + \delta t \quad \text{and} \\
 Y_{t_{0}}^{j;t_{0},x_{0};u} & = \int_{t_{0}}^{t_{0} + \delta t} \Bigl[c_j\bigl(s, X_s^{t_{0},x_{0};u}, \alpha \bigr)  + \dfrac{\partial}{\partial t} \varphi_j(s, X_{s}^{t_{0},x_{0};u}) + \mathcal{L}_{t}^{\alpha} \varphi_j(s, X_{s}^{t_{0},x_{0};u}) \notag \\
&\quad + \hat{G}_j\bigl(s, D_x \varphi_j(s, X_{s}^{t_{0},x_{0};u}) \cdot \sigma(s, X_{s}^{t_{0},x_{0};u}, \alpha)\bigr) \Bigr] d s, \quad j=1,2, \ldots, n.
\end{align*}
Further, if we substitute the above results in \eqref{EqP4.4}, we obtain
\begin{align}
\int_{t_{0}}^{t_{0} + \delta t} & \Bigl[c_j\bigl(s, X_s^{t_{0},x_{0};u}, \alpha \bigr)  + \dfrac{\partial}{\partial t} \varphi_j(s, X_{s}^{t_{0},x_{0};u}) + \mathcal{L}_{t}^{\alpha} \varphi_j(s, X_{s}^{t_{0},x_{0};u})  \notag \\
& \quad  + \hat{G}_j\bigl(s, D_x \varphi_j(s, X_{s}^{t_{0},x_{0};u}) \cdot \sigma(s, X_{s}^{t_{0},x_{0};u}, \alpha)\bigr)  \Bigr] d s  \ge 0, \quad j=1,2, \ldots, n.  \label{EqP4.6}
\end{align}
Then, dividing the above equation by $\delta t$ and letting $\delta t \rightarrow 0$, we obtain
\begin{align*}
  c_j\bigl(t_{0},x_{0}, \alpha \bigr) + \dfrac{\partial}{\partial t} \varphi_j(t_{0},x_{0}) & + \mathcal{L}_{t}^{\alpha} \varphi_j(t_{0},x_{0}) + \hat{G}_j\bigl(t_{0}, D_x \varphi_j(t_{0},x_{0}) \cdot \sigma(t_{0},x_{0}, \alpha)\bigr) \ge 0, \\
 &\hspace{1.7 in} j=1,2, \ldots, n. 
\end{align*}
Note that, since $\alpha \in W$ is arbitrary, we can rewrite the above condition as follows
\begin{align}
\dfrac{\partial}{\partial t_{0}} \varphi_j(t_{0},x_{0}) + \min_{\alpha \in U} \Bigl\{c_j\bigl(t_{0},x_{0}, \alpha \bigr) + \mathcal{L}_{t}^{\alpha} \varphi_j(t_{0},x_{0}) &+ \hat{G}_j\bigl(t_{0}, D_x \varphi_j(t_{0},x_{0}) \cdot \sigma(t_{0},x_{0}, \alpha)\bigr), \notag \\
&\hspace{1.0 in} j=1,2, \ldots, n \Bigr\} \ge 0, \label{EqP4.7} 
\end{align}
which attains its minimum in $U$ (which is a compact set in $\mathbb{R}^d)$. Thus, $V^{u}(\cdot,\cdot)$ is a viscosity subsolution of \eqref{Eq2.13}, with boundary condition $\varphi(T, x)=\Psi(T, x)$.

On the other hand, suppose that $\varphi \in C_b^{1,2}([0, T] \times \mathbb{R}^d; \mathbb{R}^n)$ and assume that $\varphi_j \le V_j^{u}$ on $[0, T] \times \mathbb{R}^d$ and $\min_{(t,x)} \bigl[V_j^{u}(t,x) - \varphi_j(t,x)\bigr]=0$ for each $j=1,2, \ldots, n$. Then, we consider a point $(t_{0},x_{0}) \in [0, T] \times \mathbb{R}^d$ so that $\varphi_j(t_{0},x_{0}) = V_j^{u}(t_{0},x_{0})$ (i.e., a local minimum at $(t_{0},x_{0})$). Further, for a small $\delta t > 0$, Let $\tilde{u}_s$, which is restricted to $\Sigma_{[t_{0},t_{0} + \delta t]}$, be an $\epsilon \delta t$-optimal control for \eqref{Eq2.17} at $(t_{0},x_{0})$. Then, proceeding in this way as \eqref{EqP4.6}, we obtain the following 
\begin{align}
\int_{t_{0}}^{t_{0} + \delta t} & \Bigl[c_j\bigl(s, X_s^{t_{0},x_{0};u}, \tilde{u}_s \bigr) + \dfrac{\partial}{\partial t} \varphi_j(s, X_{s}^{t_{0},x_{0};u}) + \mathcal{L}_{t}^{\tilde{w}_s} \varphi_j(s, X_{s}^{t_{0},x_{0};u})  \notag \\
&\quad  + \hat{G}_j\bigl(s, D_x \varphi_j(s, X_{s}^{t_{0},x_{0};u}) \cdot \sigma(s, X_{s}^{t_{0},x_{0};u}, \tilde{u}_s)\bigr) \Bigr] d s \le \epsilon \delta t, \quad j=1,2, \ldots, n.  \label{EqP4.8}
\end{align}
As a result of this, we also obtain the following
\begin{align}
\int_{t_{0}}^{t_{0} + \delta t} & \min_{\alpha \in U} \Bigl\{c_j\bigl(s, X_s^{t_{0},x_{0};u}, \alpha \bigr) + \dfrac{\partial}{\partial t} \varphi_j(s, X_{s}^{t_{0},x_{0};u}) + \mathcal{L}_{t}^{\alpha} \varphi_j(s, X_{s}^{t_{0},x_{0};u})  \notag \\
&\quad  + \hat{G}_j\bigl(s, D_x \varphi_j(s, X_{s}^{t_{0},x_{0};u}) \cdot \sigma(s, X_{s}^{t_{0},x_{0};u}, \alpha)\bigr) \Bigr\} d s \le \epsilon \delta t, \quad j=1,2, \ldots, n.  \label{EqP4.9}
\end{align}
Note that the following mappings 
\begin{align*}
(s, x, \alpha) \rightarrow  \Bigl[c_j\bigl(s, x, \alpha \bigr) + \dfrac{\partial}{\partial t} \varphi_j(t, x) + \mathcal{L}_{t}^{\alpha} \varphi_j(t, x) &+ \hat{G}_j\bigl(t, D_x \varphi_j(t, x) \cdot \sigma(t, x, \alpha)\bigr) \Bigr], \\
                                                                                        &\hspace{1.2 in} j=1,2, \ldots, n,
\end{align*}
are continuous and, since $U$ is compact, then $s \rightarrow X_s^{t_{0},x_{0};u}$ is also continuous. As a result of this, the expressions under the integral in \eqref{EqP4.9} are continuous for all $j=1,2, \ldots, n$. Further, if we divide both sides of \eqref{EqP4.9} by $\delta t$ and letting $\delta t \rightarrow 0$, then we obtain the following
\begin{align}
\dfrac{\partial}{\partial t_{0}} \varphi_j(t_{0},x_{0}) + \min_{\alpha \in U} \Bigl\{c_j\bigl(t_{0},x_{0}, \alpha \bigr) &+ \mathcal{L}_{t}^{\alpha} \varphi_j(t_{0},x_{0}) + \hat{G_j} \bigl(t_{0}, D_x \varphi_j(t_{0},x_{0}) \cdot \sigma(t_{0},x_{0},\alpha)\bigr) \notag\\
&\hspace{1.6in} j=1,2, \ldots, n \Bigr\} \le \epsilon.
 \label{EqP4.10} 
\end{align}
Notice that, since $\epsilon$ is arbitrary, we conclude that $V^{u}(\cdot,\cdot)$ is a viscosity supersolution of \eqref{Eq2.13}, with boundary condition $\varphi(T,x)=\Psi(T, x)$. This completes the proof of Proposition~\ref{P4}.
\end{IEEEproof}

\section{Further remarks} \label{S4}
In this section, we briefly comment on the implication of our main results -- when such results are implicitly used as additional information for solving optimal allocation problems. Note that, for $\xi \in L^2\bigl(\Omega, \mathcal{F}_T, \mathbb{P}, \mathbb{R}^n\bigr)$, with $\xi_j \in L^2\bigl(\Omega, \mathcal{F}_T, \mathbb{P}, \mathbb{R}\bigr)$ and $j=1,2, \ldots, n$, if we impose an additional restriction on the stochastic process $\bigl\{Y_t^{0,x;u_{\cdot}}, \,\,t \in [0, T] \bigr\}$, which is viable in $K$ for $\mathbb{P}$-{\it almost}\, $\omega \in \Omega$, to satisfy the following condition 
\begin{align*}
\mathcal{E}^{G} \bigl[\xi \cdot \mathbf{1}_n \bigl \vert \mathcal{F}_t \bigr] &= \mathcal{E}^{G} \bigl[\xi_{t,T}\bigl(u\bigr) \cdot \mathbf{1}_n \bigl \vert \mathcal{F}_t \bigr]\\
                                                                                                                &= \mathcal{E}^{G} \bigl[ Y_t^{0,x;u_{\cdot}} \cdot \mathbf{1}_n\vert \mathcal{F}_T \bigr], 
\end{align*}
for all $t \in [0, T]$ and for any admissible decision $u_{\cdot} \in \mathcal{U}_{[0,T]}$ (with $\mathbf{1}_n$ is a unit column vector with $(n \times 1)$ dimension). Then, verifying the above condition is amounted to solving an optimal allocation problem, where the target data $\xi \in L^2\bigl(\Omega, \mathcal{F}_T, \mathbb{P}, \mathbb{R}^n\bigr)$ is optimally allocated or distributed to the solution $\bigl\{Y_t^{0,x;u_{\cdot}}, \,\,t \in [0, T] \bigr\}$ of the multi-dimensional BSDE in \eqref{Eq2.6}. 

In other words, if there exists an optimal admissible decision process $u_{\cdot}^{\ast} \in \mathcal{U}_{[0,T]}$, {\it for almost all} $(s, \omega) \in [0, T] \times \Omega$, together with the corresponding solution $X_s^{t,x; u^{\ast}}$, that satisfies equation~\eqref{Eq3.15}. Then, such a solution also provides useful information to characterize all equilibrium solutions 
\begin{align*}
J\bigl[u^{\ast}\bigr] &= \xi_{0,T}(u^{\ast}) \\
                               &= \bigl[\xi_{0,T}^1(u^{\ast}), \xi_{0,T}^2(u^{\ast}), \ldots, \xi_{0,T}^n(u^{\ast})\bigr]^T, 
\end{align*}
with the following partial ordering on $K$
\begin{align*}
\mathcal{E}^{G} \bigl[\xi_{0,T}\bigl(u^{\ast}\bigr) \bigl \vert \mathcal{F}_0 \bigr] \prec \mathcal{E}^{G} \bigl[\xi_{0,T}\bigl(u\bigr) \bigl \vert \mathcal{F}_0 \bigr], \quad \forall u_{\cdot} \in \mathcal{U}_{[0,T]},
\end{align*}
when $\mathcal{E}^{G_j} \bigl[\xi_{0,T}^j\bigl(u^{\ast}\bigr) \bigl \vert \mathcal{F}_0 \bigr] \le \mathcal{E}^{G_j} \bigl[\xi_{0,T}^j\bigl(u\bigr) \bigl \vert \mathcal{F}_0 \bigr]$ for all $j=1,2, \ldots, n$, with strict inequality for at least one $j \in \{1,2, \ldots, n\}$ (e.g., see \cite{Bor62} or \cite{ChatDT00} in the context of optimal risk allocations and equilibrium solutions).

%
%
%
%
%


\begin{thebibliography}{99}

\bibitem{Pro90}
P. Protter,
{\it Stochastic integration and stochastic differential equations: a new approach},
Springer-Verlag, Berlin, Germany, 1990.

\bibitem{ParP90}
E.~Pardoux and S. Peng,
``Adapted solutions of backward stochastic differential equation,"
{\it System Control Lett.,} vol. 14, pp.~55--61, 1990.

\bibitem{Pen04} 
S. Peng,
{\it Nonlinear expectations, nonlinear evaluations and risk measures},
Lecture Notes in Mathematics. Springer, 2004.

\bibitem{CorHMP02}
F. Coquet, Y. Hu, J. M\'{e}min and S. Peng,
``Filtration-consistent nonlinear expectations and related $g$-expectations,"
{\it Probab. Theory Related Fields,} vol. 123, pp.~1--27, 2002.

\bibitem{Ros06}
E. Rosazza Gianin,
``Risk measures via $g$-expectations,"
{\it Insur. Math. Econ.}, vol. 39, pp.~19--34, 2006.

\bibitem{ArtDEH99}
P. Artzner, F. Delbaen, J. M. Eber and D. Heath,
``Coherent measures of risk,"
{\it Math. Finance}, vol. 9, pp.~203--228, 1999.

\bibitem{DetS05}
K. Detlefsen and G. Scandolo,
``Conditional and dynamic convex risk measures,"
{\it Finance Stoch.}, vol. 9, pp.~539--561, 2005.

\bibitem{Xu16}
Y. Xu,
``Multidimensional dynamic risk measure via conditional $g$-expectation,"
{\it Math. Finance,} vol. 26, pp.~638--673, 2016.

\bibitem{BuckQR00}
R. Buckdahn, M. Quincampoix and A. Rascanu,
``Viability property for a backward stochastic differential equation and applications to partial differential equations,"
{\it Prob. Theory Rel. Fields,} vol. 116, pp.~485--504, 2000.

\bibitem{AuDaPra95}
J.-P. Aubin and G. Da Prato,
``Stochastic Nagumo's viability theorem,"
{\it Stochastic Anal. Appl.,} vol. 13, pp.~1--11, 1995.

\bibitem{AuDaPra98}
J.-P. Aubin and G. Da Prato,
``The viability theorem for stochastic differential inclusions,"
 {\it Stochastic Anal. Appl.,} vol. 16, pp.~1--15, 1998.

\bibitem{Ale39} 
A. D. Alexandrov,
``The existence almost everywhere of the second differential of a convex function and some associated properties of convex surfaces,"
 {\it Ucenye Zapiski Leningrad. Gos. Univ. Ser. Math., (Russian)}, vol. 37, pp.~3--35, 1939.

 \bibitem{HUP95}
Y. Hu and S. G. Peng,
``Solutions of forward-backward stochastic differential equations,"
{\it Probab. Theory Related Fields,} vol. 103, pp.~273--283, 1995.

\bibitem{LIW14}
J. Li and Q. Wei,
``Optimal control problems of fully coupled FBSDEs and viscosity solutions of Hamilton-Jacobi-Bellman equations,"
{\it SIAM J. Control Optim.}, vol. 52, pp.~1622--1662, 2014.

\bibitem{ParT99} 
 E. Pardoux and S. J. Tang,
``Forward-backward stochastic differential equations and quasilinear parabolic PDEs,"
{\it Probab. Theory Related Fields}, vol. 114, pp. 123--150, 1999.

\bibitem{MaZZ08}
J. Ma, J. Zhang and Z. Zheng
``Weak solutions for forward-backward SDEs -- a martingale problem approach," 
{\it Ann. Probab.,} vol. 36(6, pp.~2092--2125, 2008.

\bibitem{MaPY94}
J. Ma, P. Protter, and J. M. Yong,
``Solving forward-backward stochastic differential equations explicitly - a four step scheme,"
Probab. Theory Related Fields, 98 (1994), pp.~339--359.

\bibitem{Pen92}
S. G. Peng,
``A generalized dynamic programming principle and Hamilton-Jacobi-Bellman equation,"
{\it Stoch. Stoch. Reports.}, vol. 38, pp.~119--134, 1992.

\bibitem{CraIL92}
M. G. Crandall, H. Ishii and P. L. Lions,
``User's guide to viscosity solutions of second order partial differential equations,"
{\it Bull. Amer. Math. Soc.}, vol. 27, pp.~1--67, 1992.

\bibitem{FleS06}
W. H. Fleming and H. M. Soner,
{\it Controlled Markov processes and viscosity solutions},
Springer, 2006.

\bibitem{BefVP16}
G. K. Befekadu, A. Veremyev and E. L. Pasiliao,
``Dynamic risk measures and related risk-averse decision problems,"
 in {\em Proc. 2017 American Contr. Conf.,} May 2017, Seattle, pp. 3494--3499.
 
\bibitem{DelZo95}
M. Delfour and J. P. Zolezio,
``Oriented distance functions in shape analysis and optimization," in {\it Control and optimal design of distributed parameter systems} (eds. J. E. Lagnese, \,\etal), IMA Vol. Math. Appl. vol. 70, pp.~39--72, New York, NY: Springer-Verlag, 1995. 

\bibitem{Bor62}
K. Borch,
``Equilibrium in a reinsurance market,"
{\it Econometrica}, vol. 30, pp.~424--444, 1962.

\bibitem{ChatDT00}
A. Chateauneuf, R. A. Dana and J. M. Tallon,
``Optimal risk-sharing rules and equilibria with Choquet expected-utility,"
{\it J. Math. Econ.}, vol. 34, pp.~191--214, 2000.

\end{thebibliography}
\end{document}